\documentclass[11pt]{amsart}

\usepackage[centertags]{amsmath}
\usepackage{amsfonts}
\usepackage{amssymb}
\usepackage{amsthm}
\usepackage{euscript}
\usepackage{newlfont}
\usepackage{euscript} 
\newtheorem{theorem}{Theorem}
\newtheorem{prop}{Proposition}
\newtheorem{lemma}[prop]{Lemma}

\theoremstyle{definition}
\newtheorem{definition}{Definition}

\theoremstyle{remark}
\newtheorem{remark}[theorem]{Remark}

\renewcommand{\leq}{\leqslant}
\renewcommand{\geq}{\geqslant}

\renewcommand{\cal}{\EuScript}
\newcommand{\abs}[1]{\left\vert#1\right\vert}
\newcommand{\norm}[1]{\left\Vert#1\right\Vert}
\newcommand{\R}{\mathbb R}

\newcommand{\nablash}{\nabla{\kern -.75 em
     \raise 1.5 true pt\hbox{{\bf/}}}\kern +.1 em}
\newcommand{\Deltash}{\Delta{\kern -.69 em
     \raise .2 true pt\hbox{{\bf/}}}\kern +.1 em}
\newcommand{\Rslash}{R{\kern -.60 em
     \raise 1.5 true pt\hbox{{\bf/}}}\kern +.1 em}

\newcommand{\dist}{\operatorname{dist}}

\newcommand{\D}{\partial}

\newcommand{\supp}{\operatorname{supp}}

\renewcommand{\div}{\operatorname{div}}
\newcommand{\eps}{\varepsilon}
\newcommand{\p}{\varphi}
\newcommand{\tm}{\widetilde m}
\newcommand{\tR}{\widetilde R}
\newcommand{\tA}{\tilde A}
\newcommand{\Qt}{\widetilde Q}
\newcommand{\gt}{{\tilde g}}
\renewcommand{\S}{\cal S}
\newcommand{\Mh}{\hat M}
\newcommand{\gh}{\hat g}
\newcommand{\Eh}{\hat E}
\newcommand{\Mt}{\widetilde M}
\newcommand{\Et}{\widetilde E}
\newcommand{\uh}{\hat u}
\newcommand{\N}{\cal N}
\newcommand{\Q}{\cal Q}
\newcommand{\T}{\cal T}
\newcommand{\B}{\cal B}
\newcommand{\F}{\cal F}

\newcommand{\phit}{\tilde\phi}
\newcommand{\Gammah}{\hat \Gamma}

\newcommand{\RN}{Reissner-Nordstr\"om}
\newcommand{\MP}{Majumdar-Papapetrou}

\begin{document}
\title[Penrose Inequality with Charge]{On a Penrose Inequality with Charge}

\author{Gilbert Weinstein}
\address{Department of Mathematics\\
        University of Alabama at Birmingham\\
        Birmingham, Alabama 35205}
\email{weinstein@uab.edu}
\thanks{The research of the first author
        was supported in part by
        NSF Grant~DMS-0205545.  The research of the second author
        was supported in part by
        NSF Grant~DMS-0222387.}

\author{Sumio Yamada}
\address{Mathematical Institute\\
        Tohoku University\\
        Sendai, Japan 980-8578}
\email{yamada@math.tohoku.ac.jp}

\date{\today}




\begin{abstract}
We construct a time-symmetric asymptotically flat initial data set
to the Einstein-Maxwell Equations which satisfies
\[
    m - \frac12\left(R + \frac{Q^2}{R}\right)<0,
\]
where $m$ is the total mass, $R=\sqrt{A/4\pi}$ is the area radius
of the outermost horizon and $Q$ is the total charge.  This yields
a counter-example to a natural extension of the Penrose Inequality
for charged black holes.
\end{abstract}

\maketitle

\section{Introduction}

There has recently been much interest among geometers and
mathematical relativists in inequalities bounding the total mass
of initial data sets from below in terms of other geometrical
quantities. The first such inequality is the Positive Mass
Theorem~\cite{SchoenYau79,Witten}.  We rephrase the Riemannian
version of this result as the following variational statement:
\emph{among all time-symmetric asymptotically flat initial data
sets for the Einstein-Vacuum Equations, flat Euclidean 3-space is
the unique minimizer of the total mass.}  Thus, the total mass
satisfies $m\geq0$ with equality if and only if the data set is
isometric to $\R^3$ with the flat metric.  See the next section
for precise definitions.

A stronger result is the Riemannian version of the Penrose
Inequality, which can be stated in a similar variational vein:
\emph{among all time-symmetric asymptotically flat initial data
sets for the Einstein-Vacuum Equations with an outermost minimal
surface of area $A$, the Schwarzschild slice is the unique
minimizer of the total mass.} In other words, $m\geq R/2$ where
$R=\sqrt{A/4\pi}$ is the area radius of the outermost horizon, and
equality occurs if and only if the data is isometric to the
\emph{Schwarzschild slice}:
\[
    g_{ij} = \left(1+\frac{m}{2r}\right)^4  \delta_{ij}.
\]

When these results are phrased in this fashion, a natural question
is whether similar variational characterizations of the other
known stationary solutions of the Einstein Equations hold. In
particular, one could ask whether among all asymptotically flat
axisymmetric maximal gauge initial data sets for the
Einstein-Vacuum Equations with an outermost minimal surface of
area $A$ and angular momentum $J$, the Kerr slice is the unique
minimizer of the mass.  Such a statement would imply that:
\begin{equation}    \label{eq:KerrPenrose}
    m\geq \frac12 \left(R^2 + \frac{4J^2}{R^2}\right)^{1/2}.
\end{equation}
with equality if and only if the data is isometric to the
\emph{Kerr slice}.  Since it is not known how to define the
angular momentum of a finite surface, it is necessary to assume
the axisymmetry of the data set.  With that hypothesis, if $X$ is
the generator of the axisymmetry, then the Komar integral:
\[
    J(S) = \frac{1}{8\pi} \int_S k_{ij} X^i n^j\, dA
\]
gives a quantity which depends only on the homological type of $S$
and tends to the total angular momentum, as $S$ tends to the
sphere at infinity.

A similar question can be asked with charge replacing angular
momentum: is the \RN\ slice the unique minimizer of the mass among
all asymptotically flat time-symmetric initial data sets for the
Einstein-Maxwell Equations?  This is equivalent to asking whether
the following inequality holds for any data set:
\begin{equation}    \label{eq:ChargedPenrose}
    m \geq \frac12 \left(R + \frac{Q^2}{R} \right),
\end{equation}
where $Q$ is the total charge, with equality if and only if the
data is a \RN\ slice. As above, the charge:
\[
    Q(S) = \frac1{4\pi} \int_S E_i n^i \, dA
\]
depends only on the homological type of $S$.

When the horizon is connected,
inequality~\eqref{eq:ChargedPenrose} can be proved by using the
Inverse Mean Curvature flow~\cite{HuiskenIlmanen,Jang}. Indeed,
the argument in~\cite{Jang} relies simply on Geroch montonicity of
the Hawking mass --- which still holds for the weak flow
introduced by Huisken and Ilmanen in~\cite{HuiskenIlmanen}
--- while keeping track of the scalar curvature term
$R_g=2\bigl(|E|^2+|B|^2\bigr)$. However, when the horizon has
several components the same argument yields only the following
inequality:
\[
    m \geq \frac12 \max_i \left( R_i + \frac{\left(\min \sum_i \eps_i
    Q_i\right)^2}{R_i} \right),
\]
where $R_i$ and $Q_i$ are the area radii and charges of the
components of the horizon $i=1,\dots,N$, $\eps_i=0$ or $1$, and
the minimum is taken over all possible combinations.

It is the purpose of this paper to point out
that~\eqref{eq:ChargedPenrose} does not hold. We prove:

\begin{theorem} \label{thm:main}
There is a strongly asymptotically flat time-symmetric initial
data set $(M,g,E,0)$ for the Einstein-Maxwell Equations such that:
\begin{equation}    \label{eq:ViolateChargedPenrose}
    m - \frac12\left(R + \frac{Q^2}{R}\right)<0.
\end{equation}
\end{theorem}

In 1984, Gibbons~\cite{Gibbons} conjectured an inequality similar
to~\eqref{eq:ChargedPenrose}. However, in his conjecture, the
right hand side of~\eqref{eq:ChargedPenrose} is taken to be
additive over connected components of the horizon.  Thus,
Gibbons's conjecture states that:
\begin{equation}    \label{eq:GibbonsPenrose}
    m \geq \frac12 \sum_i \left(R_i + \frac{Q_i^2}{R_i}\right).
\end{equation}
In particular, when there is no electromagnetic field this
inequality reduces to:
\begin{equation}    \label{eq:GibbonsNoCharge}
    m \geq \frac12 \sum_i R_i,
\end{equation}
which is stronger than the usual Riemannian Penrose inequality;
\[
    m\geq \frac12 \left(\sum_i R_i^2\right)^{1/2}.
\]
It is not known whether~\eqref{eq:GibbonsNoCharge} holds, but two
Schwarzschild slices a large distance apart would seem to violate
this inequality. Gibbons further conjectured that equality occurs
in~\eqref{eq:GibbonsPenrose} if and only if the data is \MP; see
the next section for a description of these metrics. We note that
these metrics do not actually have horizons and are not
asymptotically flat in the sense of Definition~\ref{def:AF}.
Instead, they have one asymptotically flat end and $N$
asymptotically cylindrical ends which we will call \emph{necks}.
The cross-sections of these necks are spheres with mean curvature
tending to zero as the surfaces goes further down the end.

Our construction is based on the fact that the \MP\ metrics
`violate'~\eqref{eq:ChargedPenrose}, say with $N=2$, and
$m_1=m_2$.  They do not strictly speaking violate~(2) since they
are not asymptotically flat and do not possess horizons.  In order
to remedy these failures, we glue two such copies along the necks.
The gluing procedure we use is an adaptation of the conformal
perturbation method developed for the vacuum case
in~\cite{IsenbergMazzeoPollack}. In fact in our setting, some of
the technical difficulties arising from the generality of the
construction in~\cite{IsenbergMazzeoPollack} are absent. However,
while it is easy to show the existence of a two-component minimal
surface in the resulting metric, we must also show
that~\eqref{eq:ChargedPenrose} is violated with $R$ the area
radius of the outermost horizon.  This requires ruling out minimal
surfaces outside the necks which we can accomplish by letting
$m\to0$ which is equivalent after rescaling to taking the two
masses in the initial \MP\ far apart.

We point out that this counter-example has little to do with the
Cosmic Censorship conjecture.  In fact, as pointed out by
Jang~\cite{Jang}, inequality~\eqref{eq:ChargedPenrose} is
equivalent to:
\[
    m - \sqrt{m^2-Q^2} \leq R \leq m + \sqrt{m^2-Q^2},
\]
and only the upper bound would follow from Cosmic Censorship using
Penrose's heuristic argument.  Our counter-example violates the
lower bound.

The paper is organized as follows.  In the next section, we define
some terms, and set-up the notation. In section~\ref{sec:gluing},
we carry out the gluing. In the last section, we show that the
parameters can be chosen so that the resulting initial data
violates~\eqref{eq:ChargedPenrose}

We wish to thank Greg Galloway, Robert Hardt, and Robert Wald for
useful discussions on this paper.  We thank the American Institute
of Mathematics for its hospitality.  The first author also thanks
the Erwin Schr\"odinger Institute for its hospitality.

\section{Preliminaries}
\label{sec:prelim}

\begin{definition}    \label{def:AF}
Let $(M,g)$ be a 3-dimensional Riemannian manifold.  We say that
$(M,g)$ is \emph{strongly asymptotically flat} (SAF) if there is a
compact set $K\subset M$ such that $M\setminus K$ is the disjoint
union of finitely many \emph{ends} $N_\nu$, $\nu=1,\dots,k$, each
end $N_\nu$ is diffeomorphic to $\R^3$ minus a ball and admits a
coordinate system in which the metric satisfies:
\[
    g_{ij} - \delta_{ij} \in C^{2,\alpha}_{-1}(N_\nu).
\]
\end{definition}

Here $C^{2,\alpha}_{-1}(N_\nu)$ denotes the class of functions
$\phi$ such that $r\abs \phi$, $r^2\abs{\D \phi}$, $r^3\abs{\D^2
\phi}$ and $r^{3+\alpha} \abs{x-y}^{-\alpha} \abs{\D^2
\phi(x)-\D^2\phi(y)}$ are bounded. While the bound is coordinate
dependent, the set of function is independent of coordinates.  We
will focus out attention on one end which we will denote by $N_+$.
We will denote all the other ends collectively as $N_-$. In fact,
in this paper we are dealing exclusively with two-ended SAF
manifolds so that $N_-$ consists of only one end. By adding a
point $\infty_-$ (or in the general case $k-1$ points) at infinity
in $N_-$ and conformally compactifying, we obtain an
asymptotically flat Riemannian manifold with one end. We now
consider the class $\S$ of smooth surfaces $S$ which bound a
compact region $\Omega$ such that $\infty_-\in\Omega$.  In this
class, it makes sense to speak of the outer unit normal.  If $S_1,
S_2\in\S$, we will say that $S_1$ \emph{encloses} $S_2$ if the
corresponding regions $\Omega_1$ and $\Omega_2$ satisfy
$\Omega_1\supset\Omega_2$.

If $(M,g)$ is strongly asymptotically flat, the \emph{total mass}
$m$ of the end $N_+$ is defined by:
\[
    m = \frac{1}{16\pi} \lim_{r\to\infty}\int_{S_r}
    \bigl(g_{ij,j}-g_{jj,i}\bigr) n^{i} \, dA,
\]
where $S_r$ is the Euclidean coordinate sphere in $N_+$, $n$ its
unit normal in $\delta$, and $dA$ the area element induced on
$S_r$ from $\delta$.

\begin{definition}
A \emph{horizon} $S$ is a minimal surface in $(M,g)$ which belongs
to $\S$. An \emph{outermost horizon} is a horizon which is not
\emph{enclosed} within any other horizon. A surface $S\in\S$ is
\emph{outer minimizing} if it has area no greater than any other
surface which encloses it.
\end{definition}

Note that for $r$ large enough, $S_r\in\S$ and has positive mean
curvature with respect to its outer unit normal.  Thus, by
minimizing area over all surfaces in $\S$ which enclose the
outermost horizon $S$, and are enclosed in $S_r$, we obtain a
minimal surface $S_1$ which encloses $S$.  It then follows from
the outermost property of $S$ that $S=S_1$; see~\cite[Theorem 1',
p.~645]{meeks-simon-yau}. We conclude that an outermost horizon is
also outer minimizing, a fact which will be used in the last
section.

A time-symmetric initial data set $(M,g,E,B)$ for the
Einstein-Maxwell Equations consists of a Riemannian manifold
$(M,g)$, and two vector fields $E$ and $B$ on $M$ such that:
\[
    R_g = 2 \bigl(\abs{E}_g^2+\abs{B}_g^2\bigr), \>
    \div_g E = \div_g B = 0, \> E \times B = 0, \>
    \int_S g(B,n_g) \, dA = 0,
\]
where $R_g$ is the scalar curvature of $g$, and $S\subset M$ is an
arbitrary closed surface with normal $n_g$ of unit length in $g$.
We say that the set $(M,g,E,B)$ is strongly asymptotically flat if
$(M,g)$ is SAF, and if $E\in C^{2,\alpha}_{-2}$, $B\in
C^{2,\alpha}_{-3}$.

Choose $N>0$, $m_k>0$, and $p_k\in\R^3$ for $k=1,\dots,N$, and let
$r_k$ denote the Euclidean distance to $p_k$ in $\R^3$. The
\emph{\MP\ }solutions are given by:
\begin{equation}
\label{eq:MP}
    u = \left(1+\sum_{k=1}^N \frac{m_k}{r_k}\right)^{1/2}, \quad
    g_{ij} = u^4\delta_{ij}, \quad
    E_i =  2 \nabla_i \log u, \quad
    B_i = 0.
\end{equation}
When $N=1$, this is simply the extreme case $m=|Q|$ of the \RN\
data set.  Note that if we take $E_-=-2\nabla\log u$ instead of
$E=2\nabla\log u$, we get another solution with charges of
opposite sign.

For simplicity, we will restrict ourselves to the case $N=2$,
$m_1=m_2=m$, i.e. $u= \bigl(1+m/r_1+m/r_2\bigr)^{1/2}$.  It is not
difficult to check that $(M,g,E,0)$ satisfies the Einstein-Maxwell
time-symmetric constraints. In fact, the metric $-u^{-4}dt^2+g$ is
a static solution of the Einstein-Maxwell equations.  Let $r$
denote the Euclidean distance from the origin. We denote by
$B_i(\rho)=\{r_i<\rho\}$ the Euclidean ball of radius $\rho$
centered at $p_i$, and by $B_0(\rho)=\{r<\rho\}$ the Euclidean
ball of radius $\rho$ centered at the origin. Note that for $R$
large enough $N=\R^3\setminus B_0(R)$ equipped with the metric $g$
is a SAF end, and the necks $B_i(\rho)\setminus\{p_i\}$ are
asymptotically cylindrical. It is easy to check that the total
mass $\mu$ of $N$ is $2m$, the total charge $Q=\int_S g(E,n)\, dA$
is $2m$, while the total cross sectional area $A$ of both necks is
asymptotically $8\pi m^2$, i.e., $R=\sqrt2 m$.  Thus, we get:
\[
    \mu-\frac12\left(R+\frac{Q^2}{R}\right)
    =2m - \frac12(\sqrt2 m + 2\sqrt 2 m) = m
    \left(2-\frac3{\sqrt2}\right) < 0.
\]
However $(M,g)$ admits no horizon.  In the next section, we remedy
this by gluing  at the necks a second copy of opposite charges.
The solution of the constraints is achieved through a conformal
perturbation argument.  We will then show in
Section~\ref{sec:outermost} that the resulting data set possesses
a horizon which violates~\eqref{eq:ChargedPenrose}.

\section{The Gluing}
\label{sec:gluing}

Let $(M_\pm,g_\pm,E_\pm,0)$ be two copies of the \MP\ data, with
$E_-=-E_+$. In this section, we show that we can glue these two
copies along their necks.  This gluing will be performed by a
perturbation method with perturbation parameter $T>0$ large.
Whenever a possible ambiguity might arise, we use a subscript (or
superscript) $+$ (or $-$ respectively) to indicate a quantity
associated with $M_+$ (or $M_-$ respectively).

For convenience, we take $p_1=(0,0,1)$ and $p_2=(0,0,-1)$. The
gluing is accomplished in three steps. In the first step, we
truncate the necks at $r_i=e^{-T}$, and introduce cut-offs in the
regions $e^{-T+1}<r^\pm_i<e^{-T+2}$ to obtain a transition to
round cylinders. This yields data on $M_+$ and $M_-$ which matches
in the regions $e^{-T}< r^\pm_i<e^{-T+1}$ of the necks. We can
then identify the corresponding boundaries $r_i^\pm=e^{-T}$ in
$M_+$ and $M_-$ creating a two-ended SAF data set
$(\Mh,\gh,\Eh,0)$. However, this data no longer satisfies the
constraint equations in the cut-off regions. In the second step,
we restore the divergence constraint $\div_{\gh}\Eh'=0$,
$\Eh'=\Eh-\nabla \p$, by solving a linear problem $\Delta_{\gh} \p
= \div_{\gh}\Eh$, $\p\to0$ at $\infty$. Finally, in the last step,
we use a perturbation argument to find a conformal deformation
$(\Mt,\gt,\Et,0)=(\Mh,\phi^4\gh,\phi^{-6}\Eh')$ which satisfies
the constraints.  It is easy to see that the divergence constraint
is automatically preserved under the above conformal
transformation $g\mapsto\phi^4g$, $\Eh'\mapsto\phi^{-6}\Eh'$,
i.e., we have $\div_\gt\Et=0$. The Gauss constraint
$R_\gt=|\Et|_\gt^2$ is then satisfied if and only if $\phi$
satisfies the following nonlinear equation~\cite{isenberg76}:
\begin{equation}    \label{eq:Lichnerowicz}
    L_{\gh} \phi = -\frac{|\Eh'|_{\gh}^2}{4\phi^3},
\end{equation}
where $L_{\gh}=\Delta_{\gh} - \frac18 R_{\gh}$ is the conformal
Laplacian of $\gh$.  Section~\ref{sec:Lichnerowicz} is therefore
devoted to showing that for $T$ large enough, there is a positive
solution $\phi$ of~\eqref{eq:Lichnerowicz} such that $\phi-1$ is
small in $C^{2,\alpha}_{-1}$.  This gluing technique is an
adaptation of ~\cite{IsenbergMazzeoPollack}.

\addtocounter{subsection}{-1}

\subsection{Function Spaces and Elliptic Theory}
\label{sec:elliptic}

Let $(M,g)$ be a SAF manifold with $K\subset M$ compact and
$M\setminus K$ the disjoint union of finitely many ends $N_\nu$.
Let $\sigma\geq1$ be a weight function on $M$ such that $\sigma=1$
on $K$, and equals the Euclidean distance $r$ on each end $N_\nu$
for $r$ large enough. Let $C^{k,\alpha}_{-\beta}(M)$ be the set of
functions $\phi$ on $M$ whose $k$-th order derivatives are
H\"older continuous and for which the norm
$\norm{\phi}_{C^{k,\alpha}_{-\beta}}$ defined below is finite:
\begin{align*}
    \norm{\phi}_{C^k_{-\beta}}&=\sum_{i=
    0}^k\norm{\sigma^{\beta+i}D^i\phi}_{C^0} \\
    [D^k \phi]_{\alpha,-\beta} &=
    \sup_{0<\dist(x,y)<\rho} \, \sigma(x,y)^{\beta+\alpha}\,
    \frac{\abs{P_y^x D^k \phi(y)- D^k \phi(x)}}{\dist(x,y)^\alpha}
    \\[1ex]
    \norm{\phi}_{C^{k,\alpha}_{-\beta}}&=\norm{\phi}_{C^k_{-\beta}}
    +[D^k\phi]_{C^{\alpha}_{-\beta-k}}.
\end{align*}
Here $D^i \phi$ represents the tensor of $i$-th order derivatives
of $\phi$, $\rho$ is the injectivity radius of $(M,g)$,
$\sigma(x,y)=\max\{\sigma(x),\sigma(y)\}$, and $P_y^x$ is parallel
translation along the shortest geodesic from $y$ to $x$.

\begin{theorem}  \label{thm:elliptic}
Let $(M,g)$ be a SAF manifold.
\begin{enumerate}
\renewcommand{\labelenumi}{(\alph{enumi})}
\item Let $\phi\in C^0_{-\beta}(M)$ and $\Delta_g \phi \in
C^{0,\alpha}_{-\beta-2}(M)$, then $\phi\in
C^{2,\alpha}_{-\beta}(M)$ and
\begin{equation}    \label{eq:CBCS:estimate}
    \norm{\phi}_{C^{2,\alpha}_{-\beta}}\leq C(\norm{\phi}_{C^0_{-\beta}}+\norm{\Delta_g
    \phi}_{C^{0,\alpha}_{-\beta-2}}).
\end{equation}
\item Let $0<\beta<1$, $\nu>2$, and let $h\in
C^{0,\alpha}_{-\nu}(M)$.  If the operator $\Delta_g-h
:C^{2,\alpha}_{-\beta}(M)\to C^{0,\alpha}_{-\beta-2}(M)$ is
injective then it is an isomorphism.
\end{enumerate}
\end{theorem}

This theorem is stated in~\cite{parkerlee}, but the reader is
referred to~\cite{choquet-bruhat79} for the proof. Unfortunately,
the proof of part~(b) in~\cite{choquet-bruhat79} has a small gap,
which is nevertheless easily remedied. For details, please refer
to~\cite[Appendix]{smithweinstein}, where a complete proof is
given for the case $M=\R^3$.  The proof for general SAF manifolds
is a straightforward combination of the arguments
in~\cite{choquet-bruhat79} and~\cite{smithweinstein}.

\subsection{Preparation}

Fix $T>0$, and let $\chi(r)$ be a smooth positive nondecreasing
cut-off function such that $\chi(r)=1$ for $r>e^{-T+2}$, and
$\chi(r)=0$ for $r<e^{-T+1}$. Let $\chi_i=\chi(r_i)$, $i=1,2$, and
define:
\[
    \uh = \left(\chi_1\chi_2 + \chi_2\,\frac{m}{r_1} +
    \chi_1\,\frac{m}{r_2}\right)^{1/2}, \quad \gh = \uh^4\delta =
    \left(\frac{\uh}{u}\right)^4 g, \quad
    \Eh_i =  2 \nabla_i \log \uh.
\]
We note that:
\begin{equation}
\label{eq:uh/u}
    \abs{1-(\uh/u)^4} \leq C
    e^{-T}, \quad \abs{\nabla\log\uh-\nabla\log u}_{g} \leq C e^{-T},
    \quad \abs{\Delta_{g} (\uh/u)} \leq C e^{-T}
\end{equation}
where $C$ is a constant independent of $T$.  This implies
\begin{equation}    \label{eq:hat-smallness}
    \abs{\gh-g}_g \leq C e^{-T}, \qquad \abs{R_{\gh} - 2
    |\Eh|^2_{\gh}} \leq C e^{-T}, \qquad \bigl|\div_{\gh}
    \Eh\bigr|\leq Ce^{-T}.
\end{equation}
Introduce the notations:
\begin{gather}
    \label{eq:defB}
    B_i(\rho)=\{r_i<\rho\},  \qquad D(\rho) = B_1(\rho)\cup B_2(\rho)\\
    \label{eq:defG}
    \Gamma_i(\rho)=\{e^{-T}\leq r_i<\rho\}, \qquad \Gamma(\rho) = \Gamma_1(\rho)\cup\Gamma_2(\rho).
\end{gather}
On $\Gamma_i(e^{-T+1})$, $\gh=m^2(dr_i^2/r_i^2 + d\omega^2)$ is a
round cylindrical metric with $d\omega^2$ the standard metric on
the unit sphere, and $\Eh=dr_i/r_i$ is parallel. Thus, if we take
two copies $M_\pm=\R^3\setminus D(e^{-T})$, then both the metrics
$\gh_\pm=\gh$ and the vector fields $\Eh_\pm=\pm\Eh$ match on
$r_i^\pm=e^{-T}$, and we can identify these boundaries to form a
doubly-connected-sum $\Mh = M_+ \# M_-$.  We will denote the
metric on $\Mh$ by $\gh$ and the vector field by $\Eh$. We note
that $(\Mh,\gh)$ is a two-ended SAF manifold.  We denote
$\Gammah_i(\rho)=\Gamma_i^+(\rho)\cup\Gamma_i^-(\rho)$ and
$\Gammah(\rho)=\Gamma^+(\rho)\cup\Gamma^-(\rho)$.  We have
suppressed the dependence on $T$ in order not to encumber the
notation.

We now fix the weight function $\sigma=\sigma(r)$ to be $1$ on
$\{r^\pm\leq 3\}$ in $M_\pm$, monotone in $r$, and equal to
$r^\pm$ on $\{r^\pm>4\}$.  In addition, we can assume that it is
even with respect to reflections across the cuts $\D
M_+=\{r_1=e^{-T}\}\cup\{r_2=r^{-T}\}$. Note that since $\uh/u=1$
outside $\Gammah(e^{-T+2})$, the quantities
in~\eqref{eq:hat-smallness} vanish outside this set, hence these
estimates hold also with any weighted norms. In particular:
\begin{equation}
\label{eq:N1}
    \norm{R_{\gh} - 2
    |\Eh|^2_{\gh}}_{C^{0,\alpha}_{-3}} \leq C e^{-T},
    \qquad
    \norm{\div_{\gh} \Eh}_{C^{0,\alpha}_{-3}} \leq C e^{-T}.
\end{equation}

Throughout the rest of this section, $C$, $C'$, $c$ will denote
various constants independent of $T$.  In order to simplify the
notation, we may at times change the value of such constants. This
abuse of notation can be justified by simply taking the maximum of
the previous and current value of the constant.

We will need the following essentially local elliptic estimate.

\begin{prop}    \label{prop:elliptic}
Let $0<\beta<1$, $\nu>2$,  and let $h\in C^{0,\alpha}_{-\nu}$
satisfy $h\geq0$. There is a constant $C$ independent of $T$, such
that for each $T$ large enough
\begin{equation}
\label{eq:local}
    \norm{\phi}_{C^{2,\alpha}_{-\beta}} \leq C\left(
    \norm{\phi}_{C^{0}_{-\beta}} + \norm{(\Delta_{\gh}-h)\phi}_{{C^{0,\alpha}_{-\beta-2}}}\right).
\end{equation}
\end{prop}

\begin{proof}
On either end $\Omega_\pm=M_\pm\setminus B_0^\pm(4)$, we can use
an argument using local estimates and the scaling of annuli as
in~\cite[Proposition 26]{smithweinstein} to get a weighted
estimate:
\[
    \norm{\phi}_{C^{2,\alpha}_{-\beta}(\Omega_\pm)} \leq C\left(
    \norm{\phi}_{C^{0}_{-\beta}(\Omega'_\pm)} +
    \norm{(\Delta_{\gh}-h)\phi}_{C^{0,\alpha}_{-\beta-2}(\Omega'_\pm)}\right).
\]
with a constant $C$ independent of $T$, where
$\Omega'_\pm=M_\pm\setminus B_0^\pm(3)$. Now, let
$K=\Mh\setminus\bigl[\{r^+\geq5\}\cup\{r^-\geq5\}\bigr]$, then $K$
can be covered by finitely many geodesic balls
$B^{\gh}_{q_i}(\rho)$ of radius $\rho>0$ sufficiently small, so
that the elliptic constant of $\gh$ written in normal coordinates
on $B^{\gh}_{q_i}(2\rho)$ is uniformly bounded above and below.
While the number of balls depends on $T$, $\rho$ can be chosen
independently of $T$.  We have local elliptic estimates:
\[
    \norm{\phi}_{C^{2,\alpha}(B^{\gh}_{q_i}(\rho))} \leq C\left(
    \norm{\phi}_{C^{0}(B^{\gh}_{q_i}(2\rho))} +
    \norm{(\Delta_{\gh}-h)\phi}_{C^{0,\alpha}(B^{\gh}_{q_i}(2\rho))}\right),
\]
where $C$ depends on $\rho$ but is independent of $i$ or $T$.
Collecting these estimates yields~\eqref{eq:local}.
\end{proof}

\subsection{The Divergence Constraint}
\label{sec:divegence}

In this section, we restore the divergence constraint by solving
the following linear problem:
\[
    \Delta_{\gh} \p = f, \qquad \p\to0 \mbox{\ at $\infty$.}
\]
where $f = \div_{\gh}(\Eh)$.  We must also ensure that $\p$ tends
to zero when $T$ tends to infinity.

\begin{prop}    \label{prop:phi}
For each $m>0$ small enough, and each $T$ large enough, there is a
unique solution $\p\in C^{2,\alpha}_{-1}$ of the equation:
\[
    \Delta_{\gh} \p = f,
\]
on $\Mh$, where $f=\div_{\gh}(\Eh)$. Furthermore,
\begin{equation}
\label{eq:vsmall}
    \norm{\p}_{C^{2,\alpha}_{-1}}\leq CT^2e^{-T}.
\end{equation}
where the constant $C$ is independent of $T$.
\end{prop}

\begin{proof}
The existence of a solution $\p\in C^{2,\alpha}_{-1}$ is standard,
see e.g.~\cite{choquet-bruhat79}. The smallness of $\p$,
inequality~\eqref{eq:vsmall}, will follow from the elliptic
estimates in Proposition~\ref{prop:elliptic} once we obtain a
weighted supremum bound:
\begin{equation}
\label{eq:weightedsup}
    \sup_{\Mh} \sigma \abs{\p} \leq CT^2e^{-T},
\end{equation}
where $C$ is independent of $T$.  This is obtained by a comparison
argument using the maximum principle.  Note that the function $f$
has $\supp f\subset \Gammah(e^{-T+2})$, and $P=\sup_T
\bigl(e^{T}\sup \abs{f}\bigr)<\infty$. Furthermore, $f$ is odd
with respect to reflection across the cuts $\D M_+$ which implies
that $\p$ is also odd, hence $\p=0$ on $\D M_+$. Now let:
\[
    \psi(r)= -e^{-T}\int_{e^{-T}}^r \frac{\log(s)}{s(s+m)}\, ds.
\]
We claim that if $m$ is small enough and $T$ is large enough, then
$w = \psi(r_1)+\psi(r_2)$ has the following properties on $M_+$:
\begin{enumerate}
\renewcommand{\labelenumi}{(\roman{enumi})}
    \item $0< w \leq m^{-1}T^2 e^{-T}$ on $M_+$.
    \item $\Delta_{\gh} w\leq 0$ on $M_+$.
    \item $\Delta_{\gh} w\leq -ce^{-T}$ on $\Gamma(e^{-T+2})$
    for some $c>0$ independent of $T$.
\end{enumerate}
These properties imply that $\p-c^{-1}P w$ satisfies
$\Delta_{\gh}(\p-c^{-1}Pw) \geq 0$, and $\p-c^{-1}Pw\leq0$ both on
$\D M_+$ and at $\infty$.  Thus, we get $\p\leq c^{-1}Pw$ on
$M_+$. Similarly, by considering the function $\p+c^{-1}Pw$, we
get $\p\geq-c^{-1}Pw$ on $M_+$. This yields an unweighted supremum
estimate:
\begin{equation}
\label{eq:sup}
    \sup_{M_+} \abs{\p}\leq \frac{P}{cm}T^2e^{-T}.
\end{equation}
By symmetry, the same estimate holds on $M_-$. Now, in order to
get the weighted estimate~\eqref{eq:weightedsup}, let
$\Omega=M_+\setminus B^+_0(3)$, and let $v$ be the solution of the
following problem:
\[
    \Delta_{g} v=0 \mbox{\ in $\Omega$}, \qquad v=1 \mbox{\ on $\D\Omega$},\qquad
    v\to0 \mbox{\ at $\infty$}.
\]
There is a constant $C$ such that $0<v\leq C\sigma^{-1}$. Let
$K=P/cm$, then the functions $\pm\p+KT^2e^{-T}v$ are harmonic in
$\Omega$ with respect to $g=\gh$, are non-negative on $\D\Omega$,
and tends to $0$ at $\infty$, hence by the maximum principle
$\pm\p+KT^2e^{-T}v\geq0$ in $\Omega$. Hence, we obtain
$\sigma\abs{\p}\leq KCT^2e^{-T}$ on $\Omega$. Combining
with~\eqref{eq:sup}, the weighted estimate~\eqref{eq:weightedsup}
follows.

It remains to prove the claims (i)-(iii).  Denote
$\psi_i=\psi(r_i)$, and note that
\[
    \max\psi_i=\psi(1)=-e^{-T}\int_{e^{-T}}^1\frac{\log s}{s(s+m)}\,
    ds \leq \frac1{2m}\,T^2e^{-T},
\]
whence $w\leq m^{-1}T^2e^{-T}$.  A similar estimate shows that
$\psi(1)>T^2e^{-T}/4m$.  On the other hand
\[
    \psi(1)-\psi(\infty) = e^{-T}\int_1^\infty \frac{\log s}{s(s+m)}\, ds \leq
    e^{-T} < \psi(1),
\]
provided $T$ is large enough.  We conclude that $w>0$ if $T$ is
large enough proving~(i).  In order to establish~(ii) and~(iii),
we first note that it is sufficient to prove these with $\gh$
replaced by $g$.  Indeed, suppose that~(i) and~(ii) hold with $g$
instead of $\gh$.  Then we have:
\[
    \Delta_{\gh} w = \left(\frac u{\uh}\right)^4 \left(
    \Delta_g w + \left(\frac u{\uh}\right)^2
    g\left(\nabla(\uh/u)^2,\nabla w\right) \right)
\]
and $\uh/u=1$ outside $\Gamma(e^{-T+2})$, while on
$\Gamma(e^{-T+2})$:
\[
    \abs{\nabla w}_g \leq \frac{C}{m^2}\, Te^{-T},
    \qquad
    \abs{\nabla(\uh/u)^2}_g \leq C e^{-T}.
\]
It follows that~(i) and~(ii) also hold with $\gh$ once we replace
$c$ by say $c/2$, provided $T$ is large enough.

We now turn to proving ~(i) and~(ii) with respect to $g$.  Let
$g_1=u_1^4\delta$ be the one-black-hole \MP\ metric, i.e.,
$u_1^2=1+m/r_1$.  One easily calculates:
\[
    \Delta_{g_1}\psi_1 = -\frac{e^{-T}}{(r_1+m)^3}.
\]
Let $\theta=\nabla r_1\cdot\nabla r_2$ denote the inner product of
$\nabla r_1$ and $\nabla r_2$ with respect to $\delta$, then:
\begin{align*}
    \Delta_g\psi_1 &= u^{-6}\,\div_\delta(u^2\nabla\psi_1) \\
        &= u^{-6}\,
        \div_\delta\left(\left(u_1^2+
        \frac{m}{r_2}\right)\nabla\psi_1\right) \\
        &= \left(\frac{u_1}{u}\right)^6\, \Delta_{g_1}\psi_1 +
        \frac{m}{u^6}\,\div_\delta\left(\frac{1}{r_2u_1^2}\,u_1^2\nabla\psi_1\right)\\
        &=-\left(\frac{u_1}{u}\right)^6\left(1+\frac{m}{r_2u_1^2}\right)\frac{e^{-T}}{(r_1+m)^3}
        + \frac{mu_1^2}{u^6}\,
        \nabla\left(\frac1{r_2u_1^2}\right)\cdot\nabla\psi_1 \\
        &=-\left(\frac{u_1}{u}\right)^6\left(1+\frac{m}{r_2u_1^2} +
        \frac{m^2\log r_1}{r_2(r_1+m)u_1^2} - \frac{mr_1^2\log
        r_1}{r_2^2(r_1+m)}\,\theta\right) \,
        \frac{e^{-T}}{(r_1+m)^3}.
\end{align*}
Note that $r_1>1$ on $B_2(1)$ hence, there, we can estimate
$\Delta_{g}\psi_1$ above by the only positive term on the
right-hand side:
\[
        \Delta_{g}\psi_1\leq \abs{\frac{me^{-T}r_2\log
        r_1}{(r_2+mr_2/r_1+m)^3r_1(r_1+m)}\,\theta}  \leq
        \frac{e^{-T}\abs{\log r_1}}{m^2r_1^2}\, r_2.
\]
Furthermore, as $r_1\to0$, then $\Delta_{g}\psi_1\to-e^{-T}/m^3$.
It follows that, provided $m<1$, one can choose $\eps>0$
independent of $T$ and $m$ such that:
\[
    \Delta_{g}\psi_1 \leq
    \begin{cases}
        e^{-T}/4m^3, & \text{when $r_2<\eps$}\\
        -e^{-T}/2m^3, &\text{when $r_1<\eps$}.
    \end{cases}
\]
Now, we can choose $m>0$ small enough, so that
$\Delta_{\gh}\psi_1\leq0$ when $r_1, r_2\geq\eps$.  By symmetry,
we have analogous estimates for $\Delta_{g}\psi_2$.  We conclude
that:
\begin{gather*}
    \Delta_{g} w = \Delta_g \psi_1 + \Delta_g \psi_2 \leq 0,
    \qquad \text{when $r_1,r_2\geq\eps$.} \\
    \Delta_g w \leq -\frac{e^{-T}}{4m^3}, \qquad \text{on
    $\Gamma(\eps)$.}
\end{gather*}
Properties~(ii) and (iii) now follow provided $T>-\log\eps+2$.
This completes the proof of Proposition~\ref{prop:phi}.
\end{proof}

Defining $\Eh'=\Eh-\nabla \p$, we now have $\div_{\gh}(\Eh')=0$,
and in view of~\eqref{eq:vsmall} and~\eqref{eq:N1}:
\begin{equation}
\label{eq:N1'}
    \norm{R_{\gh} - 2
    |\Eh'|^2_{\gh}}_{C^{0,\alpha}_{-3}} \leq C T^2e^{-T}.
\end{equation}

\subsection{The Gauss Constraint}
\label{sec:Lichnerowicz}

In this section, we prove that for each $T$ large enough, there is
a positive solution $\phi\in 1+C^{2,\alpha}_{-1}(M)$ of
Equation~\eqref{eq:Lichnerowicz}. We first prove the following
estimate which gives a uniform bound on the inverse of the
linearized operator associated with~\eqref{eq:Lichnerowicz}.  The
proof is adapted from~\cite[Proposition 8]{IsenbergMazzeoPollack}.

\begin{prop}    \label{prop:InvN}
Let
\[
    h = \frac18 R_{\gh} + \frac34 |\Eh'|_{\gh}^2,
\]
Then there is a constant $C$ independent of $T$ such that if
$\phi\in C^{2,\alpha}_{-2/3}$, then
\[
    \norm{\phi}_{C^{2,\alpha}_{-2/3}} \leq C
    \norm{(\Delta_{\gh}-h)\phi}_{C^{0,\alpha}_{-8/3}}.
\]
\end{prop}

\begin{remark} For $T$ bounded, this follows from Theorem~2.
\end{remark}

\begin{proof}
Suppose the contrary.  Then there is a sequence $T_j\to\infty$,
and $\phi_j\in C^{2,\alpha}_{-2/3}$ satisfying
\begin{equation}
\label{eq:contradiction}
    \norm{\phi_j}_{C^{2, \alpha}_{-2/3}} = 1 \quad \forall j, \qquad
    \norm{(\Delta_{\gh}-h)\phi_j}_{C^{0, \alpha}_{-8/3}}  \to 0 \mbox{\ as $j \to
    \infty$}.
\end{equation}
By Proposition~1, we have
\[
    \norm{\phi_j}_{C^{2, \alpha}_{-2/3}} \leq C \big(
    \norm{\phi_j}_{C^{0}_{-2/3}} + \norm{(\Delta_{\gh}-h)\phi_j}_{C^{0,
    \alpha}_{-8/3}} \big),
\]
with $C$ independent of $j$.  Hence, in view
of~\eqref{eq:contradiction}, we have $\eps
> 0$ such that
\begin{equation}
\label{eq:epsilon}
    \eps \leq \norm{\phi_j}_{C^{0}_{-2/3}} \leq 1,
\end{equation}
for all $j$.   We now consider the following two cases:
\begin{itemize}
\item[(i)] There is $\tau>0$ such that for any $\delta>0$, we
have:
\[
    \limsup_j\norm{\phi_j}_{C^0(\Gammah(\delta))}
    \geq \tau,
\]
\item[(ii)] For every $\tau>0$, there exists $\delta>0$ so that:
\[
    \limsup_j\norm{\phi_j}_{C^0(\Gammah(\delta))}<\tau.
\]
\end{itemize}
Note that $\Gammah(\delta)$ is the union of the two necks cut at
$r_i^\pm=\delta$.

\emph{Case}~(i). For each integer $k$ large enough, take
$\delta_k=e^{-k}$ in~(i).  Then there is $j_k$ large enough so
that $T_{j_k}>k$ and so that there exists
$p_k\in\Gammah(\delta_k)$ with $\abs{\phi(p_{j_k})}\geq\tau/2$.
Without loss of generality we may assume that
$p_{j_k}\in\Gamma_1^+(\delta_k)$. Furthermore, by passing to a
subsequence, we may assume that $j_k=k$, i.e., $T_k>k$, and
$p_k\in\Gamma_1^+(\delta_k)$.  We define a coordinate $s$ on $\Mh$
by:
\[
    s = \pm (\log r_1 + T) \mbox{ on } M_{\pm}.
\]
Denote $s_k=s(p_k)$, then it follows that $0\leq s_{k}<T_{k}-k$.
Now let
\[
    \Lambda_k = \{ s_{k} - k/2 < s < s_{k} + k/2 \}.
\]
The part $\{ s = s_{k} + k/2 \}$ of the boundary of $\Lambda_k$
has $r^+_1$ coordinate equal to:
\[
    \exp({s_{k} + k/2 - T_{k}}) <  e^{-k/2} \to 0, \quad \text{as
    $k\to\infty$}.
\]
A similar estimate holds for the other part of the boundary.  It
follows that the metric $\gh$ on $\Lambda_k$ converges to
$m^2(ds^2 + d\omega^2)$ as $k\to\infty$.  Let $(\Lambda,g_0)$
denote the standard round cylinder with the metric $g_0=m^2(ds^2 +
d\omega^2)$. We will identify the points on $\Lambda_{k}$ with
those of $\Lambda$ via the identity map induced by the $(s,
\omega)$ coordinates. Observe that $\cup\Lambda_k=\Lambda$. Using
the compactness of the embedding $C^{2,
\alpha'}(\Lambda_k)\hookrightarrow C^{2,\alpha}(\Lambda_k)$, $0
\leq \alpha' < \alpha$, we can now select a subsequence, which we
now denote $\phi_k$ again, and a function $\phi_0$ on $\Lambda$
such that $\phi_k\to\phi_0$ in $C^{2,\alpha'}(\Lambda_k)$ for each
fixed $k$. Furthermore, there is a point $p_0$ in the
cross-section $\{s=0\}$ of $\Lambda$ such that
$\abs{\phi_0(p_0)}\geq\tau/2$, hence $\phi_0$ is not identically
zero. The scalar curvature $R_{\gh}$ on $\Lambda_k$ converges to
$2/m^2$ and $|\hat{E}'|_{\gh}^2$ converges to $1/m^2$. Thus, the
coefficients of $\T=\Delta_{\gh}-h$ converge uniformly on compact
sets to the coefficients of
\[
    \T_0 = \Delta_{g_0} - \frac{1}{m^2}.
\]
Hence we get $\T \phi_{k} \to \T_0 \phi_0$ in $C^{0}(\Lambda_k)$.
Since we also have $\T\phi_k\to0$ in $C^{0,\alpha}(\Lambda_k)$
by~\eqref{eq:contradiction}, we conclude that $\phi_0$ satisfies
the linear equation
\[
\Delta_{g_0} \phi_0 - \frac{1}{m^2} \phi_0 = 0
\]
on $\Lambda$.  Since $\phi_0$ is nontrivial, it has exponential
growth in $s$ either as $s\to\infty$ or as $s\to-\infty$, in
contradiction to~\eqref{eq:epsilon}.

\emph{Case}~(ii).  We begin the treatment of this case with the
following lemma:

\begin{lemma}   \label{lemma:A}
Suppose $\phi_j$ satisfies~\eqref{eq:contradiction}
and~\emph{(ii)}, and let $A^+_\delta\subset M_+$ be the twice
perforated ball $A_\delta=B_0(3)\setminus D(\delta/2)$. Then for
each $\delta>0$, there holds $\norm{\phi_j}_{C^{1,
\alpha}(A_\delta)}\to0$ as $j\to\infty$.
\end{lemma}

\begin{proof}
Suppose not, and let $j_k$ be a subsequence such that $\phi_{j_k}$
converges to $\phi_0$ in $C^{1,\alpha'}(A_\delta)$,
$\alpha'<\alpha$.  Then $\phi_0$ is not identically zero on
$A_\delta$, hence since $h>0$ on $A_\delta$, we have:
\[
    \lim_k \int_{A_\delta}
    h\phi_{j_k}^2
    = \int_{A_\delta} h\phi_0^2 > 0.
\]
We now proceed to show that
\begin{equation}
\label{eq:limsup}
    \limsup_j \int_{A_\delta}
    h\phi_{j}^2 = 0,
\end{equation}
leading to a contradiction.  Without loss of generality, we may
assume, by passing to a subsequence, that:
\[
    \int_{A_\delta}
    h\phi_{j}^2 \to
    \limsup_j \int_{A_\delta}
    h\phi_{j}^2.
\]
If $\chi$ is any smooth cut-off function of compact support in
$M_+$, with $0\leq\chi\leq1$, and $\chi=1$ in $A_\delta$, then:
\begin{align*}
    \int_{A_\delta}
    h\phi_{j}^2
    & \leq \int_{M_+}
    \chi^2\left(\abs{\nabla\phi_{j}}_{\gh}^2 +
    h\phi_{j}^2\right)\\
    & =  -\int_{M_+}\chi^2\,\phi_j\,(\Delta_{\gh}-h)\phi_j -
    2 \int_{M_+} \chi\,\phi_j\, \gh(\nabla\chi,\nabla\phi_j).
\end{align*}
We will now choose cut-offs $\chi_k$, and a subsequence $j_k$
along which both of the terms on the right-hand side tend to zero,
proving~\eqref{eq:limsup}. By~(ii), for each integer $k$ we can
choose $0<\delta_k<\delta$ such that $\abs{\phi_j}<1/k$ on
$\Gammah(\delta_k)$ for all $j$ large enough.  Now, we can choose
$\chi_k$ supported on $B_0(k)\setminus D(\delta_k/2)$, with:
\[
    \supp\nabla\chi_k\subset \bigl[B_0(k)\setminus B_0(3)\bigr] \cup
    \bigl[D(\delta_k)\setminus D(\delta_k/2)\bigr],
\]
and satisfying:
\begin{gather*}
    \abs{\nabla\chi_k}_{\gh} \leq C/k, \qquad \text{on $B_0(k)\setminus B_0(3)$},\\
    \abs{\nabla\chi_k}_{\gh}\leq C, \qquad \text{on $D(\delta_k)\setminus
    D(\delta_k/2)$}.
\end{gather*}
Finally, by~\eqref{eq:contradiction}, we can choose $j_k>j_{k-1}$,
so that
\[
    \left(\int_{\R^3\setminus B_0(3)}\sigma^{-10/3}\right)
    \norm{(\Delta_{\gh}-h)\phi_{j_k}}_{C^{0,\alpha}_{-8/3}}\leq
    \frac1k.
\]
It then follows that:
\[
    - \int_{M_+}\chi_k^2\,\phi_{j_k}\,(\Delta_{\gh}-h)\phi_{j_k}
    < \frac1k,
\]
and:
\begin{align*}
    - \int_{M_+} \chi_k\phi_{j_k}\gh(\nabla\chi_k,\nabla\phi_{j_k})
    &\leq C \int\limits_{D(\delta_k)\setminus D(\delta_k/2)} \abs{\phi_{j_k}}
    + Ck^{-1} \int\limits_{B_0(k)\setminus B_0(3)} \sigma^{-7/3}\\
    &\leq Ck^{-1} + Ck^{-1}\int_3^k \frac{dr}{r^{1/3}}\\
    &\leq
    C\left(k^{-1}+ k^{-1/3}\right).
\end{align*}
This completes the proof of Lemma~\ref{lemma:A}.
\end{proof}

Now choose $\delta>0$ so that
\begin{equation}
\label{eq:eps}
    \limsup_j\norm{\phi_j}_{C^0(\Gammah(\delta))}<\eps,
\end{equation}
where $\eps$ is defined by~\eqref{eq:epsilon}, and define a new
manifold $(M_*, g_*)$ diffeomorphic to ${\R}^3$ by extending
smoothly the metric $\gh$ on $\R^3\setminus D(\delta/2)$ over
$D(\delta/2)$. Then extend smoothly to $D(\delta/2)$ also the
potential function $h$ so that the extended potential $h_*$
satisfies $h_*\geq0$ on $\R^3$. Let $\chi$ be a smooth cut-off
function on $\R^3$ with $0\leq\chi\leq1$, $\chi=1$ outside
$D(\delta)$, and $\chi=0$ on $D(\delta/2)$. Taking the values of
$\phi_j$ from $M_+$, we can view $\chi\phi_j$ as a function on
$M_*$, and we find:
\[
   (\Delta_{g_*}-h_*)\chi\phi_j
    = (\Delta_{\gh} - h)\chi \phi_j
    = \chi (\Delta - h)\phi_j +
    2 \gh(\nabla\chi,\nabla \phi_j) + \phi_j \Delta_{\gh}\chi.
\]
Hence, we can estimate:
\[
    \norm{(\Delta_{g_*}-h_*)\chi\phi_j}_{C^{0,\alpha}_{-8/3}(M_*)}
    \leq \norm{(\Delta - h)\phi_j}_{C^{0,\alpha}_{-8/3}} +
    C\norm{\phi_j}_{C^{1,\alpha}(A_\delta)} \to 0,
\]
by~\eqref{eq:contradiction} and Lemma~\ref{lemma:A}. It then
follows by Theorem~2 part~(b) applied to $(M_*,g_*)$ that:
\[
    \norm{\chi\phi_j}_{C^{2,\alpha}_{-2/3}(M_*)} \to 0.
\]
Thus, we obtain:
\[
    \norm{\phi_j}_{C^{0}_{-2/3}
    (M_+\setminus\Gamma^+(\delta))} \to 0.
\]
Similarly, we obtain:
\[
    \norm{\phi_j}_{C^{0}_{-2/3}
    (M_-\setminus\Gamma^-(\delta))} \to 0,
\]
and it follows that:
\[
    \norm{\phi_j}_{C^{0}_{-2/3}
    (\Mh\setminus\Gammah(\delta))} \to 0.
\]
Combining with~\eqref{eq:eps}, we conclude that:
\[
    \limsup_j \norm{\phi_j}_{C^{0}_{-2/3}
    (\Mh)} < \eps,
\]
in contradiction to~\eqref{eq:epsilon}.  This completes the proof
of Proposition~\ref{prop:InvN}.
\end{proof}

We can now prove the main result of this section.

\begin{prop} \label{thm:Lichnerowicz}
For each $m>0$ small enough and for each $T$ large enough there is
a solution $\phi\in 1+C^{2,\alpha}_{-1}(M)$ of:
\begin{equation}    \label{eq:Gauss}
    L_{\gh} \phi = -\frac{|\Eh'|_{\gh}^2}{4\phi^3} .
\end{equation}
Furthermore, as $T\to\infty$, this solution satisfies
$\norm{\phi-1}_{C^{2,\alpha}_{-1}}\to0$.
\end{prop}

\begin{proof}
Let $\N\colon 1+C^{2,\alpha}_{-2/3}\to C^{0,\alpha}_{-8/3}$ be the
following nonlinear operator:
\begin{equation}
\label{eq:defN}
    \N(1+\psi) = L_{\gh}(1+\psi) + \frac{|\Eh'|_{\gh}^2}{4 (1+\psi)^{3}}.
\end{equation}
The linearization of $\N$ about $1$ is:
\[
    d\N = L_{\gh} - \frac34 |\Eh'|_{\gh}^2 \colon C^{2,\alpha}_{-2/3}\to
    C^{0,\alpha}_{-8/3},
\]
and according to Proposition~\ref{prop:InvN}, $d\N^{-1}\colon
C^{0,\alpha}_{-8/3}\to C^{2,\alpha}_{-2/3}$ is uniformly bounded,
i.e., there is a constant $C$ independent of $T$ such that
\[
    \norm{d\N^{-1}\psi}_{C^{2,\alpha}_{-2/3}}\leq
    C\norm{\psi}_{C^{0,\alpha}_{-8/3}}.
\]
Now consider the `quadratic part' of $\N$:
\[
    \Q(\psi) = \N(1+\psi)- \N(1) - d\N(\psi)
\]
We have:
\[
    \Q(\psi) =
    \frac{|\Eh'|_{\gh}^2(6+8\psi+3\psi^2)}{4(1+\psi)^3}\,
    \psi^2,
\]
hence it follows that there is a constant $C$ independent of $T$
such that if $\eta>0$ is sufficiently small, and
$\norm{\psi}_{C^{2,\alpha}_{-2/3}}<\eta$, then the following
holds:
\begin{align}
    \label{eq:into}
    \norm{\Q(\psi)}_{C^{0,\alpha}_{-8/3}} &\leq C
    \eta^2\\
    \label{eq:contract}
    \norm{\Q(\psi_1)-\Q(\psi_2)}_{C^{0,\alpha}_{-8/3}} &\leq 2C\eta
    \norm{\psi_1-\psi_2}_{C^{2,\alpha}_{-2/3}}.
\end{align}
Now, choose $0<\lambda<1$, $\eta>0$ such that $\eta<\lambda/2C^2$,
and $T>0$ such that $T^2e^{-T}<\eta^2$.  Then, if $\B$ is the ball
of radius $\eta$ in $C^{2,\alpha}_{-2/3}$,  the map $\F$ given by:
\[
    \F(\psi) = - d\N^{-1}\bigl( \N(1) + \Q(\psi)\bigr)
\]
maps $\B$ into $\B$ and is a contraction.  Indeed, in view
of~\eqref{eq:N1'} and~\eqref{eq:into}, we have:
\[
    \norm{\F(\psi)}_{C^{2,\alpha}_{-2/3}} \leq C\left(\norm{\N(1)}_{C^{0,\alpha}_{-8/3}} +
    \norm{\Q(\psi)}_{C^{0,\alpha}_{-8/3}}\right)
    \leq C^2 (T^2e^{-T} + \eta^2) < \eta,
\]
and in view of~\eqref{eq:contract}
\begin{multline*}
    \norm{\F(\psi_1)-\F(\psi_2)}_{C^{2,\alpha}_{-2/3}} \leq
    C\norm{\Q(\psi_1)-\Q(\psi_2)}_{C^{0,\alpha}_{-8/3}} \\
    \leq 2C^2\eta
    \norm{\psi_1-\psi_2}_{C^{2,\alpha}_{-2/3}} <\lambda\norm{\psi_1-\psi_2}_{C^{2,\alpha}_{-2/3}}.
\end{multline*}
It follows that $\F$ has a fixed point $\psi$ in $\B$ which
satisfies
\[
    \N(1+\psi)=\N(1)+d\N(\psi)+\Q(\psi)=0.
\]
Furthermore, note that if $T\to\infty$, one can choose $\eta\to0$.
Thus we have $\norm{\psi}_{C^{2,\alpha}_{-2/3}}\to0.$  It also
follows from~\eqref{eq:defN} that:
\[
    \Delta_{\gh}\psi = \frac18 R_{\gh}\psi + \frac{|\Eh'|_{\gh}^2}{4(1+\psi)^3} - \frac18
    R_{\gh}.
\]
We will now show that the right-hand side above tends to zero as
$T\to\infty$ in $C^{0,\alpha}_{-3}\cap L^1$.  Indeed, we have:
\begin{gather*}
    \frac18 R_{\gh}\psi + \frac{|\Eh'|_{\gh}^2}{4(1+\psi)^3} - \frac18
    R_{\gh} \\
    = \frac18 R_{\gh}\psi - \frac18\left[ R_{\gh} -
    2|\Eh'|_{\gh}^2\right] -
    \frac{|\Eh'|_{\gh}^2}{4(1+\psi^3)}\left((1+\psi)^3-1\right)
    \\
    = \frac18 R_{\gh}\psi - \frac18\left[ R_{\gh} -
    2|\Eh|_{\gh}^2\right] - \frac14 \gh(2\Eh+\nabla\p,\nabla\p) +
    \frac{|\Eh'|_{\gh}^2(3+3\psi+\psi^2)}{4(1+\psi^3)}\,\psi
\end{gather*}
We now proceed to check that each of the terms above tends to zero
in $C^{0,\alpha}_{-3}\cap L^1$ as $T\to\infty$.  The second term
above tends to zero by~\eqref{eq:N1}, and the fact that it is
supported on a set of uniformly bounded volume. For the other
three terms, we use the fact, that if $f_i\in
C^{0,\alpha}_{-\beta_i}$, $i=1,2$, and $\beta_1+\beta_2>3$, then
\[
    \norm{f_1f_2}_{C^{0,\alpha}_{-3}\cap L^1} =
    \norm{f_1f_2}_{C^{0,\alpha}_{-3}} + \norm{f_1f_2}_{L^1} \leq
    C\norm{f_1}_{C^{0,\alpha}_{-\beta_1}}\norm{f_2}_{C^{0,\alpha}_{-\beta_2}}.
\]
If one of the factors on the right-hand side of the inequality
tends to zero and the other is bounded, then the left-hand side of
the inequality tends to zero.  The first and last term above are
of the form $f\psi$, with $\norm{f}_{C^{0,\alpha}_{-4}}$ bounded
and $\norm{\psi}_{C^{0,\alpha}_{-2/3}}\to0$.  The third term is of
the form $f|\nabla\p|_{\gh}$ with $\norm{f}_{C^{0,\alpha}_{-2}}$
bounded and $\norm{\nabla\p}_{C^{0,\alpha}_{-2}}\to0$.  We
conclude that:
\begin{equation}
\label{eq:rhs->0}
    \norm{\Delta_{\gh}\psi}_{C^{0,\alpha}_{-3}\cap L^1}\to0
\end{equation}

The result will now follow from the following Lemma.

\begin{lemma}   \label{lemma:1decay}
Suppose that $\psi\in C^{0,\alpha}_{-2/3}$ and
$\Delta_{\gh}\psi\in C^{0,\alpha}_{-3}\cap L^1$.  Then there is a
constant $C$ independent of $T$ such that:
\[
    \norm{\psi}_{C^{2,\alpha}_{-1}} \leq C\left(
    \norm{\Delta_{\gh}\psi}_{C^{0,\alpha}_{-3}\cap L^1}
    + \norm{\psi}_{C^{2,\alpha}_{-2/3}}\right).
\]
\end{lemma}

\begin{proof}[Proof of Lemma~\ref{lemma:1decay}]
The proof of this lemma is based on the proof of Proposition~29
in~\cite[Appendix]{smithweinstein}.  There, it is proved that if
$v$ is a function on $\R^3$ with $v\in C^{0,\alpha}_{-2/3}$, and
$\Delta_{\gh} v\in C^{0,\alpha}_{-3}\cap L^1$, then:
\begin{equation}    \label{eq:Prop29}
    \norm{v}_{C^{2,\alpha}_{-1}} \leq C \norm{\Delta_{\gh}
    v}_{C^{0,\alpha}_{-3}\cap L^1}.
\end{equation}
Let $\Omega_R=\{r>R\}\subset M_+$.  Clearly, for any finite
$R\geq3$, the two norms
$\norm{\cdot}_{C^{2,\alpha}_{-1}(M_+\setminus\Omega_R)}$ and
$\norm{\cdot}_{C^{2,\alpha}_{-2/3}(M_+\setminus\Omega_R)}$ are
equivalent.  Let $\chi$ be a smooth cut-off function with
$0\leq\chi\leq1$, $\chi=0$ on $M\setminus\Omega_3$ and $\chi=1$ on
$\Omega_4$.  Then $v=\chi\psi$ can be viewed as a function on
$\R^3$, and $v\in C^{2,\alpha}_{-2/3}$.  We have:
\[
    \Delta_{\gh}v = \chi\Delta_{\gh}\psi +
    2\gh(\nabla\chi,\nabla\psi) + \psi\Delta_{\gh}\chi.
\]
The last two terms above are supported on the annulus
$\Omega_3\setminus\Omega_4$, hence we can estimate:
\[
    \norm{2\gh(\nabla\chi,\nabla\psi) +
    \psi\Delta_{\gh}\chi}_{C^{0,\alpha}_{-3}\cap L^1} \leq
    C\norm{\psi}_{C^{2,\alpha}_{-2/3}},
\]
while for the first term we clearly have:
\[
    \norm{\chi\Delta_{\gh}\psi}_{C^{0,\alpha}_{-3}\cap L^1} \leq
    C\norm{\Delta_{\gh}\psi}_{C^{0,\alpha}_{-3}\cap L^1}.
\]
Thus, we obtain:
\[
    \norm{\Delta_{\gh}v}_{C^{0,\alpha}_{-3}\cap L^1} \leq C
    \left( \norm{\Delta_{\gh}\psi}_{C^{0,\alpha}_{-3}\cap L^1}
    + \norm{\psi}_{C^{2,\alpha}_{-2/3}}\right).
\]
We now conclude from~\eqref{eq:Prop29} that
\[
    \norm{\psi}_{C^{2,\alpha}_{-1}(\Omega_4)} \leq
    C\norm{\Delta_{\gh}v}_{C^{0,\alpha}_{-3}\cap L^1}
    \leq C\left( \norm{\Delta_{\gh}\psi}_{C^{0,\alpha}_{-3}\cap L^1}
    + \norm{\psi}_{C^{2,\alpha}_{-2/3}}\right).
\]
Hence, we have:
\begin{align*}
    \norm{\psi}_{C^{2,\alpha}_{-1}(M_+)} &\leq C
    \left(\norm{\psi}_{C^{2,\alpha}_{-1}(\Omega_4)} +
    \norm{\psi}_{C^{2,\alpha}_{-1}(M_+\setminus\Omega_5)}\right)\\  &\leq
    C\left( \norm{\Delta_{\gh}\psi}_{C^{0,\alpha}_{-3}\cap L^1}
    + \norm{\psi}_{C^{2,\alpha}_{-2/3}}\right)
\end{align*}
A similar estimate holds on $M_-$.  This completes the proof of
Lemma~\ref{lemma:1decay}.
\end{proof}

Taking $\phi=1+\psi$, we see that $\phi$ satisfies
Equation~\eqref{eq:Gauss}, and by~\eqref{eq:rhs->0} and
Lemma~\ref{lemma:1decay} we have
$\norm{\phi-1}_{C^{2,\alpha}_{-1}}\to0$ as $T\to\infty$. This
completes the proof of Proposition~\ref{thm:Lichnerowicz}.
\end{proof}

We have shown that for each $m>0$ sufficiently small and for each
$T$ sufficiently large, there is a two-ended solution
$(\Mt,\gt,\Et,0)$ of the Einstein-Maxwell constraints, with
$\Mt=\Mh$, $\gt = \phi^4 \uh^4 \delta = \phit^4 g$, where
$\phit=\phi\uh/u$ and $g$ is the \MP\ solution.  We note that for
any $\eta>0$, we can assure that:
\begin{equation}
\label{eq:small}
    \Vert\phit-1\Vert_{C^{2,\alpha}_{-1}} < \eta, \qquad
    \Vert E-\Et\Vert_{C^{1,\alpha}_{-2}} < \eta
\end{equation}
by taking $T$ large enough. For the sake of simplicity, we now
rename $\phit$ to be $\phi$. Furthermore, we note that this
solution admits an involutive, charge-reversing, symmetry with
fixed-point set $\Sigma_0$.

\section{The Outermost Horizon}
\label{sec:outermost}

In this section, we show that with $m>0$ fixed and sufficiently
small, we can adjust the perturbation parameter $\eta>0$ to be
small enough so that the area $\tA$ of the outermost horizon in
the conformal perturbation $\gt$ is no greater than $ 8\pi
\lambda^2m^2$ where $\lambda-1>0$ is arbitrarily small, i.e.,
$\tR\leq\sqrt2\lambda m$. Furthermore, if the perturbation
parameter $\eta>0$ is small enough we can assure that the total
mass $\tm$ of $\gt$ is no greater than $2\lambda m$ and that the
charge $\Qt$ of $\Et$ satisfies $\Qt\geq Q/\lambda$, where $Q=2m$.
Now, the function $f_Q(x)=x+Q^2/\lambda^2x$ is non-increasing for
$0<x<Q/\lambda$. Thus, if we choose $\lambda$ so that
\[
    1<\lambda<\left(\sqrt{2}-1/2\right)^{-1/4}<2^{1/4},
\]
we get:
\[
    \tm - \frac12 \left(\tR + \frac{\Qt^2}{\tR}\right) \leq
    \tm - \frac12 \left(\tR + \frac{Q^2}{\lambda^2\tR}\right) \leq
    \lambda m\left(2 - \frac1{\sqrt2} - \frac{\sqrt2}{\lambda^4}\right)
    < 0,
\]
This proves Theorem~\ref{thm:main}.

We will use the following elementary lemma.

\begin{lemma}   \label{lemma:vboundsH}
Let $M$ be a Riemannian manifold, and let $S\subset M$ be a
compact hypersurface with unit normal $n$. Let $v$ be a smooth
function on $M$, and suppose that the maximum of $v$ over $S$ is
taken at a point $q\notin\D S$ where $\nabla v\ne0$.  Then
\begin{equation}    \label{eq:lowerH}
    \abs{H(q)} \geq \left. \frac{\Delta v - \nabla^2_n v}{\abs{\nabla v}}\right|_q,
\end{equation}
where $H$ is the mean curvature of $S$.
\end{lemma}

\begin{proof}
Let $\Deltash$ denote the Laplacian with respect to the metric
induced on $S$.  We have at $q$:
\[
    \Delta v = \nabla^2_n v + H \nabla_n v + \Deltash v
    \leq \nabla^2_n v + H \nabla_n v,
\]
since $\Deltash v \leq0$ there, and since, without loss of
generality, we may take $n=\nabla v/\abs{\nabla v}$.  Thus, we
obtain:
\[
    \Delta v \leq \nabla^2_n v + H \abs{\nabla v},
\]
and~\eqref{eq:lowerH} follows,
\end{proof}

The right hand side of~\eqref{eq:lowerH} is easily recognized as
the mean curvature of the level set of $v$ at $q$.  Thus this
lemma is simply another version of the familiar geometric fact
that when two surfaces are tangent at $q$ and one lies entirely on
one side from the other, then their mean curvature at $q$ are
correspondingly ordered.  We prefer the statement in the lemma
since it simplifies some of the explicit computations below.  A
first application is the following lemma.

\begin{lemma}   \label{lemma:HnotinB3}
Let $p_0=(0,0,0), p_{1}=(0,0,1), p_2=(0,0,-1) \in \R^3$, and let
$0<\eps<1/3$. Then, for any compact surface $S\subset
B_0(3)\setminus D(\eps)$ such that $\D S\subset \D D(\eps)$, there
holds:
\[
    \sup_{S} \abs{H} \geq \frac1{6},
\]
where $H$ is the mean curvature of $S$.
\end{lemma}

\begin{proof}
We consider two cases: (i) $\max_S v>0$; and~(ii) $v\leq0$ on $S$.
For case~(i), take $v=x^2+y^2-z^2/2$. Let $q\in S$ be such that
$v(q)=\max_S v$.  Since $v<0$ on $\D D(\eps)$ we conclude that
$q\notin \D S$. Since $\Delta v = 3$, $\nabla^2_nv \leq 2$, and
$\abs{\nabla v} \leq 2r<6$, Lemma~\ref{lemma:vboundsH} now yields:
\[
    H(q) \geq \frac16.
\]
Now in case~(ii) note that since $S$ is smooth, it is contained in
the double cone $v<0$. Without loss of generality $S_1=S \cap
\{z>0\}\ne\emptyset$, and we now take $v=r_1^2=x^2+y^2+(z-1)^2$,
and let $q$ as above be such that $v(q)=\max_S v$. If $q\in\D
S_1\subset \D B_1(\eps)$, then $S_1\subset\D B_1(\eps)$ and
$H=2/\eps\geq6$ at every interior point of $S$.  On the other
hand, if $q\notin\D S_1$, then since $\Delta w = 6$, $\nabla^2_n w
= 2$, and $\abs{\nabla w} = 2r_1$, we obtain from
Lemma~\ref{lemma:vboundsH} that:
\[
    H(q) \geq \frac{2}{r_1} > \frac1{2}.
\]
\end{proof}

\begin{prop}    \label{prop:inr1}
If $m$ is sufficiently small, then for each $\eps>0$ there is
$\eta>0$ such that if $\norm{\phi-1}_{C^{2,\alpha}_{-1}}<\eta$,
then any closed surface $S\subset M$ which is minimal in the
conformal perturbation $(M,\phi^4g,\phi^{-6}E)$ of $(M,g,E)$ is
contained in $D(\eps)$.
\end{prop}

\begin{proof}
The proof is established in three stages.  We first show that $S$
cannot enter the region outside $B_0(3)$.  We then do the same for
for the twice-perforated ball $B_0(3)\setminus D(1/4)$. Finally,
we prove the result in each of the two balls $B_1(1/4)$ and
$B_2(1/4)$ separately. We will use the Euclidean metric $\delta$,
the \MP\ metric $g=u^4\delta$, and also its perturbation
$\gt=\phi^4 g$. In order to avoid confusion we will use the dot
product to denote the inner product with respect to $\delta$, and
indicate other metric objects by subscripts.  We denote $\nu=\phi
u$, and note that
\begin{equation}    \label{eq:conformalH}
    H_\gt = \div_\gt(n_\gt) = \frac1{\nu^6} \div_\delta(\nu^4 n_\delta) =
    \frac1{\nu^2} H_\delta + \frac4{\nu} \gt (\nabla\nu, n_\gt),
\end{equation}
where $H_\gt$ and $H_\delta$ denote the mean curvatures of $S$ in
the metrics $\gt$ and $\delta$ respectively.

Suppose first that $\max_S r\geq3$ where $r$ is the Euclidean
distance from $p_0$.  Then, in view of $\Delta_\delta r=2/r$,
$\abs{\nabla r}_\delta=1$, $\nabla^2 r=0$, we have according to
Lemma~\ref{lemma:vboundsH} that at the point $q$ with maximum $r$:
\[
    \abs{H_\delta(q)} \geq \frac2r.
\]
Now, $u^2=1+m/r_1+m/r_2$, hence outside $B_0(3)$, we have:
\[
    \abs{\nabla \log u}_{g} = \frac12 \frac{\abs{\nabla
    u^2}_\delta}{u^4}\leq
    \frac12 \frac{m(1/r_1^2+1/r_2^2)}{(1+m/r_1+m/r_2)^2} \leq
    \frac{m}{(r-1)^2} \leq \frac{3m}{4r}.
\]
Thus, using~\eqref{eq:conformalH} and $\abs{\phi-1},
r^2\abs{\nabla \phi}_g \leq
\norm{\phi-1}_{C^{2,\alpha}_{-1}}<\eta$, we can estimate:
\begin{align*}
    \abs{H_{\gt}(q)}
&   \geq \frac{\abs{H_\delta}}{\nu^2} -
    4\abs{\nabla\log \nu}_\gt \\
&   \geq\frac{2}{r(1+2m/3)(1+\eta)^2} -
    \frac{4}{\phi^2} \left(\abs{\nabla\log\phi}_{g}
    +\abs{\nabla\log u}_{g} \right) \\
&   \geq \frac{2}{r(1+2m/3)(1+\eta)^2} -
    \frac{4}{(1-\eta)^3}
    \left(\frac{\eta}{r^2} + \frac{3m}{4r} \right) \\
&   \geq \frac2r \left(\frac{1}{(1+2m/3)(1+\eta)^2} -
    \frac{4\eta + 3m}{2(1-\eta)^3} \right)
\end{align*}
Clearly if $m$ and $\eta$ are small enough, then
$\abs{H_{\gt}(q)}>0$, a contradiction.  We conclude that $S\subset
B_0(3)$.

Suppose now that $S$ enters $B_0(3)\setminus D(1/4)$.  Then a
similar estimate yields a point $q$ in that region where:
\[
    \abs{H_\delta(q)} \leq 4\left(1+\frac{2m}{3}\right)
    \left(\frac{\eta}{9(1-\eta)} + 16m \right).
\]
Hence, if $\eta$ and $m$ are small enough, then we have
$\abs{H_\delta(q)}<1/6$ in contradiction to
Lemma~\ref{lemma:HnotinB3}.

Therefore, we can now fix $m$ and $\eta_0$ small enough such that
if $\eta<\eta_0$ then $S$ must lie in $D(1/4)$. Consider the
closed surface $S_1=S\cap B_1(1/4)$ with $H_\gt=0$. As above, we
can estimate:
\[
    \abs{H_g} \leq 4 \abs{\nabla\log\phi}_g <  \frac{4\eta}{1-\eta}.
\]
We will now apply Lemma~\ref{lemma:vboundsH} to the surface $S_1$
and the function $r_1$ in $B_1(1/4)$ equipped with the metric $g$.
Let $q$ be the point where $r_1(q)=\max_{S_1} r_1<1/4$. We compute
at $q$, using $n_g=\nabla r_1/\abs{\nabla r_1}_g$:
\begin{gather*}
    \Delta_g r_1 = \frac1{u^6} \div_\delta(u^2\nabla r_1) =
    \frac1{u^6} \left( \frac{2u^2}{r_1} + \nabla u^2 \cdot \nabla
    r_1 \right) \\
    \nabla^2_{n_g} r_1 = \nabla_{n_g}\abs{\nabla r_1} = \nabla_{n_g} u^{-2} =
    - \frac{g(\nabla u^2, \nabla r_1)}{u^4 \abs{\nabla r_1}} =
    -\frac{\nabla u^2 \cdot \nabla r_1}{u^6}.
\end{gather*}
Thus, using $m<2$, we can estimate:
\begin{align*}
    \frac{\nabla r_1 - \nabla_{n_g}^2 r_1}{\abs{\nabla r_1}}
&   = \frac1{u^4} \left(\frac{2u^2}{r_1} + 2 \nabla u^2 \cdot
      \nabla r_1
    \right) \\
&   = \frac1{u^4} \left( \frac{2u^2}{r_1} - 2m \left(
    \frac1{r_1^2} + \frac{\nabla r_1\cdot\nabla r_2}{r_2^2} \right)\right)  \\
&   \geq \frac1{u^4} \left( \frac2{r_1} - 2m \right) \\
&   \geq \frac{r_1}{(1+m)^2}.
\end{align*}
We now obtain from Lemma~\ref{lemma:vboundsH}:
\[
    \frac{r_1(q)}{(1+m)^2} \leq \abs{H_g(q)} \leq
    \frac{4\eta}{1-\eta}.
\]
Therefore, with $m$ fixed, we see that $\max_{S_1}r_1\to0$ as
$\eta\to0$.  The same argument can be applied to $S_2=S\cap
B_2(1/4)$.  This proves Proposition~\ref{prop:inr1}
\end{proof}

\begin{prop}
Let $\lambda>1$.  Then for $m$ and $\eta$ sufficiently small, the
mass $\tm$ of $\gt$, the area $\tA$ of the outermost horizon in
$\gt$, and the charge $\Qt$ satisfy:
\begin{equation}    \label{eq:lambda}
    \tm \leq 2 \lambda m, \qquad \tA \leq 8 \pi \lambda^2 m^2,
    \qquad \Qt\geq Q/\lambda.
\end{equation}
\end{prop}

\begin{proof}
The metric of the perturbed space is $\gt=\phi^4 u^4 \delta$,
therefore its mass $\tm$ is:
\[
    \tm = 2m - \frac1{2\pi} \lim_{r\to\infty} \int_{S_r}
    \frac{\D\phi}{\D r}\, dA_0,
\]
where $A_0$ is the area element of $S_r$ in the flat metric
$\delta$.  Since $\D\phi/\D r$ can be estimated on $S_r$ by
$r^{-2} u^2 \norm{\phi}_{C^{2,\alpha}_{-1}} < r^{-2} u^2 \eta$, we
find:
\[
    \tm \leq 2m \left( 1 + \frac{\eta}{m}\right).
\]
Thus, $\tm\leq2\lambda m$ provided $\eta\leq m(\lambda-1)$.
Similarly $\Qt\geq Q/\lambda$ follows from~\eqref{eq:small}.

We now note that $\gt$ admits one horizon.  Indeed, the surface
$\Sigma_0=\{s=0\}$ cutting the neck at its midpoint is totally
geodesic, since it is the fixed-point set of the isometry sending
any point $p$ on one side of it to the corresponding point on the
other.  In particular, this surface is minimal and encloses the
end $\infty_-$, hence it is a horizon.

Now let $S$ be the outermost horizon.  Then $S$ is outer
minimizing.  According to Proposition~\ref{prop:inr1}, if $m>0$
and $\eta>0$ are sufficiently small, then $S\subset D(\eps)$,
where $\eps\leq 4\eta(1+m)^2/(1-\eta)$. Thus, $\D D(\eps)$
encloses $S$, and we conclude:
\begin{multline*}
    \tA = A_\gt(S) \leq A_\gt \bigl(\D D(\eps)\bigr) =
    \int_{\D D(\eps)} \phi^2 u^4\, dA_0 \\
    \leq (1+\eta)^4 \int_{\D D(\eps)} u^4
    dA_0 \leq 8\pi m^2 (1+\eta)^4 \left(1 +
    \eps\left(1+\frac{1}{m}\right) \right)^2 \\
    \leq 8\pi m^2 (1+\eta)^4 \left(1 +
    \frac{4\eta(1+m)^3}{m(1-\eta)}
    \right)^2.
\end{multline*}
Thus, with $m>0$ fixed and small enough, we can choose $\eta>0$
small enough to satisfy~\eqref{eq:lambda}.
\end{proof}

\bibliographystyle{amsplain}

\end{document}